# A PRIORCONDITIONED LSQR ALGORITHM FOR LINEAR ILL-POSED PROBLEMS WITH EDGE-PRESERVING REGULARIZATION

SIMON R. ARRIDGE[*], MARTA M. BETCKE[†], AND LAURI HARHANEN[‡]

**Abstract.** This article presents a method for solving large-scale linear inverse problems regularized with a nonlinear, edge-preserving penalty term such as the total variation or Perona–Malik. In the proposed scheme, the nonlinearity is handled with lagged diffusivity fixed point iteration which involves solving a large-scale linear least squares problem in each iteration. Because the convergence of Krylov methods for problems with discontinuities is notoriously slow, we propose to accelerate it by means of priorconditioning. Priorconditioning is a technique which embeds the information contained in the prior (expressed as a regularizer in Bayesian framework) directly into the forward operator and hence into the solution space. We derive a factorization-free priorconditioned LSQR algorithm, allowing implicit application of the preconditioner through efficient schemes such as multigrid. The resulting method is matrix-free i.e. the forward map can be defined through its action on a vector. We demonstrate the effectiveness of the proposed scheme on a three-dimensional problem in fluorescence diffuse optical tomography using algebraic multigrid preconditioner.

**1. Introduction.** Inverse problems arise in almost all fields of science when details of a model have to be determined from a set of observed data. Formally, we consider a mapping $A : X \to Y$ as the *forward problem* and the inversion of this mapping as the *inverse problem*. A defining characteristic of such problems is that they are *ill-posed* whenever the forward mapping is a compact operator which, for infinite-dimensional linear operators, is equivalent to a gradual decay to zero for the singular values of the mapping $A$. As a consequence of the ill-posedness of $A$, its inversion is unstable and requires *regularization*. For further details, see e.g. [3, 14, 25].

In this work, we are interested in linear problems where $X = \mathbb{R}^{n_X}$, $Y = \mathbb{R}^{n_Y}$, $n_X, n_Y \in \mathbb{N}$. We consider $X$ as a space of images (after appropriate discretization), which, for spatial dimensions greater than one, means that the problems are *large-scale*. Examples of such problems include image deblurring, denoising and inpainting, tomography based on the Radon Transform (XCT) or Attenuated Radon Transform (PET, SPECT), or fluorescence Diffuse Optical Tomography (fDOT). In a discrete implementation, the representation of $A$ would entail an explicit matrix construction that becomes infeasible even for moderately sized images due to memory restrictions. Therefore in this paper we focus on methods which use an implicit representation of $A$, which is sometimes referred to as a *matrix-free* approach. Examples include projection operators in XCT, PET and SPECT, and the solution of direct and adjoint partial differential equations in fDOT.

**1.1. Bayesian framework.** In this paper, we make use of the Bayesian framework which recasts the original reconstruction problem into a probabilistic setting by considering the unknown parameter $f$ as a random variable, cf. [9, 27, 47]. The main goal in Bayesian inference is to combine findings obtained from measurements with the knowledge that was available before any observations. Technically, such information is encoded as two probability distributions: the *likelihood* $\pi_Y(g|f)$ and the *prior* $\pi_X(f)$. The former describes the uncertainty in measurements for fixed, known parameter values, and the latter contains the information on $f$ that is independent of the measurements. An application of Bayes' formula specifies the *posterior* distribution

$$(1.1) \qquad \pi(f|g) \propto \pi_Y(g|f)\pi_X(f),$$

which quantifies all the available information on $f$.

Assuming a zero-mean Gaussian noise model, the unnormalized likelihood is given as

$$\pi_Y(g|f) \propto \exp\left\{-\frac{1}{2}\|g - Af\|^2\right\}.$$

[*]Department of Computer Science, University College London (`Simon.Arridge@cs.ucl.ac.uk`)
[†]Department of Computer Science, University College London (`M.Betcke@cs.ucl.ac.uk`)
[‡]Department of Mathematics and Systems Analysis, Aalto University (`lauri.harhanen@aalto.fi`)

Without loss of generality, we have assumed that the noise covariance has been incorporated in $A$ and $g$.

For image reconstruction problems, the prior $\pi_X$ models a distribution of images and a typical assumption is that edges should be enhanced. Such priors are intimately related to nonlinear regularization functionals which are frequently applied through solution of the corresponding Euler–Lagrange equations of nonlinear diffusion type. We therefore consider discretized versions of regularizers having the form

$$\mathcal{R}(f) := \int_\Omega r(|\nabla f(x)|) \mathrm{d}x \tag{1.2}$$

and the related priors

$$\pi_X(f) \propto \exp\{-\tau R(f)\},$$

where $R$ represents a discretization of $\mathcal{R}$, $r \in C^1([0, \infty))$ and $\tau > 0$ controls the strength of the prior. Three well known examples are

- first order differential Tikhonov (spatially homogeneous smoothing prior)

$$r(t) = \frac{1}{2}t^2, \tag{1.3}$$

- total variation [40] (sparsity prior for edges) and its smoothed approximation

$$r(t) = t, \quad \text{or} \quad r(t) = T\sqrt{1 + (t/T)^2}, \tag{1.4}$$

- Perona–Malik [38] (edge-preserving prior)

$$r(t) = \frac{1}{2}T^2 \log(1 + (t/T)^2), \quad \text{or} \quad r(t) = \frac{1}{2}T^2(1 - \exp(-t^2/T^2)). \tag{1.5}$$

Note that the penalty term corresponding to total variation (1.4) is not differentiable without the smoothing imposed by the parameter $T > 0$. For (1.5), $T$ is a threshold parameter indicating a level of image structure below which edges are considered as noise; we apply the same interpretation to this parameter in (1.4) to illustrate the generic approach, even though in the total variation literature this factor is usually stated as a purely numeric correction.

Since exploration of the posterior distribution (1.1) often entails computationally highly expensive Monte Carlo simulations, it is common to use simpler estimators. Here, we use the *Maximum A-Posteriori* (MAP) estimate

$$f_{\mathrm{MAP}} = \arg\max_{f \in X} \pi(f|g),$$

which can also be obtained by solving

$$\min_{f \in X} \left[ \Phi(f) := \frac{1}{2}\|g - Af\|^2 + \tau R(f) \right], \tag{1.6}$$

since $\Phi$ is the negative logarithm of the (unnormalized) posterior. Due to the aforementioned connection, the terms prior and regularization are used as synonyms throughout the paper.

**1.2. Krylov methods for linear inverse problems.** We first consider the case where the prior is also a zero-mean, Gaussian distribution with a (possibly improper) covariance defined via its inverse $C_X^{-1} = \tau L^\mathrm{T} L$. This corresponds to $R(f) = \frac{1}{2}\|Lf\|^2$ and simplifies (1.6) into

$$\min_{f \in X} \left[ \Phi(f) := \frac{1}{2}\|g - Af\|^2 + \tau \frac{1}{2}\|Lf\|^2 \right], \tag{1.7}$$

whence the MAP estimate for the inverse problem is equivalent to the solution of the minimization problem

$$\min_{f \in X} \left\| \begin{bmatrix} A \\ \sqrt{\tau} L \end{bmatrix} f - \begin{bmatrix} g \\ 0 \end{bmatrix} \right\|. \tag{1.8}$$

A generally applicable method for problems of this type is LSQR [37, 36] which minimizes the functional $\Phi$ over the Krylov space

$$\text{span}\left\{ A^T g, \left(A^T A + \tau L^T L\right) A^T g, \ldots, \left(A^T A + \tau L^T L\right)^{i_{\max} - 1} A^T g \right\} \subset X. \tag{1.9}$$

We note the following difficulties with methods based on the space (1.9):
- If the regularization functional is designed to promote edges, solutions exhibit slow convergence due to the slow build up of high frequencies. This is discussed in detail in Section 2.
- Regularization can be controlled either through the parameter $\tau$ or the limit on the dimension $i_{\max}$ of the solution subspace in (1.9). The former is part of the Bayesian formalism, and the latter results from purely numerical considerations of ill-posed matrix inversion. While the truncation index $i_{\max}$ is usually determined implicitly within an iterative algorithm by the choice of a stopping rule, change of $\tau$ requires a recomputation of the subspace (1.9) from scratch.

**1.3. Data and solution priorconditioning.** Priorconditioning [6, 7, 8, 9] is a technique rooted in the Bayesian framework. Assuming $L$ is invertible, the change of basis

$$\hat{f} = Lf, \quad \hat{A} = AL^{-1},$$

transforms (1.7) into

$$\min_{\hat{f} \in LX} \left[ \hat{\Phi}(\hat{f}) := \frac{1}{2} \|g - \hat{A}\hat{f}\|^2 + \tau \frac{1}{2} \|\hat{f}\|^2 \right]. \tag{1.10}$$

Since the transformed variable $\hat{f}$ has now standard multivariate normal distribution, it is sometimes referred to as *whitened* in the signal processing literature, cf. [9] for a Bayesian view on whitening.

From the algebraic perspective, (1.10) can be thought of as a transformation to the standard form. On structured meshes, it is natural to define the prior by constructing the matrix $L$ directly. In this case $L$ is typically nonsquare and hence noninvertible. Thus numerical methods, cf. [12, 28], based on the transformation to standard form need to employ the A-weighted pseudoinverse [12]. In contrast, we aim to solve problems where the prior is given through a scaled inverse of the prior covariance, $M = L^T L$, instead of the factor $L$ needed by the pseudoinverse-based methods. This is a typical situation for problems on unstructured grids.

For the transformed problem (1.10), LSQR seeks the solution from the transformed Krylov space

$$\text{span}\left\{ M^{-1} A^T g, \left(M^{-1} A^T A\right) M^{-1} A^T g, \ldots, \left(M^{-1} A^T A\right)^{i_{\max} - 1} M^{-1} A^T g \right\}. \tag{1.11}$$

Application of LSQR to (1.10) instead of (1.7) has a number of benefits:
- Local structure in the prior is embodied directly in the transformed operator $\hat{A}$, allowing quick convergence.
- The corresponding Krylov space is the same for all values of $\tau$ due to translation invariance.
- As the prior is embodied into the transformed operator, it is possible to exploit the prior even without explicit regularization, i.e., setting $\tau = 0$. The regularization can be achieved implicitly through early stopping of a Krylov solver.

- The transformed variables are *dimensionless*, which removes issues with physical units if further transformations, e.g. exponentiation, are used.

These advantages of the priorconditioned formulation motivate the current work. Taking a naïve approach, the construction of the priorconditioned problem (1.10) requires an explicit symmetric factorization $L^\mathrm{T} L = M$ with an invertible $L$, e.g., the Cholesky decomposition. For large-scale problems such a factorization may be too expensive to obtain, store and apply: Although $M$ is sparse, $L$ will have considerable fill-in which may result in a prohibitively high computational cost. Moreover, for nonlinear problems the factorization would be computed at each linearization step of a nonlinear iteration, further aggravating the computational issues. On the other hand, $M := (\tau C_X)^{-1} = L^\mathrm{T} L$ is typically a second order differential operator (see Section 1.4) which can be cheaply stored due to its sparse representation. In this work we propose a factorization-free priorconditioned LSQR algorithm, which solves the priorconditioned formulation (1.10) without actually factorizing $M$.

**1.4. Nonlinear Regularization.** The use of a nonlinear regularizer such as (1.4)–(1.5) leads to a nonlinear optimization problem (1.6), which we solve by searching for a critical point. As we aim to work with unstructured grids, it is natural to define the gradient of $R$ as a discretization of the Fréchet derivative of $\mathcal{R}$. The Fréchet derivative of $\mathcal{R}$ at $f$ in the direction $h$ is defined through

$$(1.12) \qquad \mathcal{R}'(f)h = \int_\Omega \frac{r'(|\nabla f(x)|)}{|\nabla f(x)|} \nabla f(x) \cdot \nabla h(x) \mathrm{d}x,$$

which clearly is the weak form of the inhomogeneous diffusion operator

$$(1.13) \qquad -\nabla \cdot c_f(x) \nabla$$

acting on $f$, where $c_f(x) = \frac{r'(|\nabla f(x)|)}{|\nabla f(x)|}$. We note here that when images are defined on regular structured (i.e. tensor) grids, then regularizers of type (1.2) lead to discretization schemes based on simple finite difference operators, cf. [51], that are naturally given in the factorized form (with noninvertible factors). On unstructured grids this situation is reversed: The discretization of (1.12) yields the regularization matrix $M$ in an unfactorized form and factorizing it would be prohibitively expensive as noted in Section 1.3.

Regularizers of the type (1.2) are nonlinear through their dependence on the solution. As the problems considered in this work are linear, the regularizers are the only nonlinear part of (1.6). This suggests a successive linearization approach: at each nonlinear iteration step, the regularizer is evaluated at the currently available approximation of the solution. The here presented algorithm derives from a method of this type, the lagged diffusivity fixed point iteration, introduced in [52] originally for total variation image denoising.

**1.5. Overview of the contribution.** We present a *matrix and factorization-free* algorithm for *large-scale* linear inverse problems with *nonlinear regularizers*. The resulting nonlinearity is limited to the regularization term, and is handled with lagged diffusivity fixed point iteration [52]. The solution of the intermediate linear problems is accelerated with factorization-free priorconditioning of LSQR, which allows for using efficient algorithms such as multi-grid methods for applying the preconditioner. In particular, the proposed technique avoids the costly computation of the factorization for the regularizer $M = L^\mathrm{T} L$.

The remainder of this paper is organized as follows. In Section 2, we review the idea of priorconditioning of Krylov spaces. Section 3 derives the factorization-free priorconditioned LSQR algorithm for solving the transformed least squares problem (1.10). Section 4 discusses the nonlinear solver based on successive linearizations and subspace priorconditioning for solution of the nonlinear regularized problem (1.6). The performance of the method is tested in Section 5 on a three-dimensional fDOT problem. We conclude with a summary of the results and a discussion of prospective research directions.

**2. Priorconditioning of Krylov spaces.** In this section, we review priorconditioning [6, 7, 8, 9] for accelerating the convergence of Krylov subspace methods for solving (1.7). To simplify the exposition, we revert to the corresponding normal equation

$$(A^{\mathrm{T}}A + \tau M) f = A^{\mathrm{T}}g, \tag{2.1}$$

but the numerical method introduced in Section 3 works directly on the least squares problem (1.7). Although the formulation (1.7) uses the factorization $M = L^{\mathrm{T}}L$, the algorithm introduced in Section 3 does not require such.

**2.1. Symmetric priorconditioning.** For regularizers of the form (1.2), the matrix $M$ (cf. (1.13)) and the corresponding normal equation (2.1) are symmetric, which implies that it is beneficial to use a symmetric solver such as the conjugate gradient method (CG). The convergence of CG is governed by the distribution of the eigenvalues of the system matrix, with the eigenvectors corresponding to the extremal eigenvalues converging the fastest in the associated Krylov subspace. As an example, the operator considered in Section 2.2 is a convolution whose eigenvectors are Fourier-type modes with the frequency growing as the corresponding eigenvalue decreases. Such a modal basis provides only a poor approximation for the edges even after many iterations, which is the reason why Krylov methods converge slowly for the deconvolution of functions with jumps.

For well-conditioned problems, the aim of preconditioning is to cluster the eigenvalues of the transformed problem in order to improve the convergence of Krylov subspace methods. For ill-posed inverse problems, the mapping $A$ is a discretization of a compact operator whose singular values accumulate at zero, rendering such clustering impossible. Therefore, some different form of preconditioning is necessary for inverse problems. One such approach, proposed in [23, 26], is based on an idea of splitting the solution space into two subspaces: a small subspace for representing the approximate regularized solution and its large complement. The small projected problem is solved using a direct algorithm while the remaining solution component is computed with an iterative method. Another method, proposed in [39], solves the least squares problem (1.7) by projecting it into a generalized Krylov subspace spanned by the vectors $\{g, Ag, Lg, A^2g, ALg, LAg, L^2g, \ldots\}$. While this technique incorporates the prior information into the solution space and is thus closely connected to priorconditioning, its motivation is distinctly different. Finally, we mention earlier work on preconditioning with linear differential operators [21, 24], cf. the nonlinear inhomogeneous diffusion operator (1.13) used in this work for priorconditioning.

This article follows the philosophy of priorconditioning [6, 7, 8, 9], where the preconditioner is based on the properties of the expected solution rather than those of the forward mapping. Priorconditioning provides a way of *incorporating the jump information directly into the forward operator*, which would otherwise require a large number of iterations to be built up in the Krylov subspace. As a result, a priorconditioned Krylov iteration converges significantly faster than an unpriorconditioned one.

To illustrate this behavior, we construct three Krylov spaces, which correspond to three different variants of the normal equation (2.1):

- **Unregularized:** regularization arises solely from early stopping,

$$A^{\mathrm{T}}Af = A^{\mathrm{T}}g. \tag{2.2}$$

  The solution space spanned by the Krylov solver is (cf. (1.9))

  $$\mathcal{K}^{A^{\mathrm{T}}A} = \mathrm{span}\left\{A^{\mathrm{T}}g, A^{\mathrm{T}}AA^{\mathrm{T}}g, \ldots, (A^{\mathrm{T}}A)^{i_{\max}-1}A^{\mathrm{T}}g\right\}.$$

- **Regularized (but unpriorconditioned):** regularization through the choice of $\tau$ (and possibly early stopping),

  $$(A^{\mathrm{T}}A + \tau M) f = A^{\mathrm{T}}g.$$

The solution space spanned by the Krylov solver is

$$\mathcal{K}^{A^{\mathrm{T}}A+\tau M} = \mathrm{span}\left\{A^{\mathrm{T}}g, (A^{\mathrm{T}}A+\tau M)A^{\mathrm{T}}g, \ldots, (A^{\mathrm{T}}A+\tau M)^{i_{\max}-1}A^{\mathrm{T}}g\right\}.$$

- **Priorconditioned:** symmetrically (split) preconditioned normal equation

$$(2.3) \qquad (L^{-\mathrm{T}}A^{\mathrm{T}}AL^{-1} + \tau I)\hat{f} = L^{-\mathrm{T}}A^{\mathrm{T}}g, \qquad f = L^{-1}\hat{f},$$

where we have used a factorization $M = L^{\mathrm{T}}L$ with an invertible factor $L$. The solution space spanned by the Krylov solver is (cf. (1.11))

$$\mathcal{K}^{M^{-1}A^{\mathrm{T}}A} = \mathrm{span}\left\{M^{-1}A^{\mathrm{T}}g, M^{-1}A^{\mathrm{T}}AM^{-1}A^{\mathrm{T}}g, \ldots, \right.$$
$$\left. (M^{-1}A^{\mathrm{T}}A)^{i_{\max}-1}M^{-1}A^{\mathrm{T}}g\right\}.$$

Note that for the priorconditioned formulation, the solution $f$ is contained in $\mathcal{K}^{M^{-1}A^{\mathrm{T}}A}$ while the actual subspace built by the Krylov solver for the transformed solution $\hat{f}$ is $L\mathcal{K}^{M^{-1}A^{\mathrm{T}}A}$. Due to the shift invariance of the Krylov spaces, $\mathcal{K}^{M^{-1}A^{\mathrm{T}}A} = \mathcal{K}^{M^{-1}A^{\mathrm{T}}A+\tau I}$ for any choice of $\tau$. Hence, in the priorconditioned least squares problem (2.3) $\tau$ acts solely as a damping parameter.

It is illuminating to look at the Lanczos vectors (orthonormalized Krylov basis vectors) in each of the above spaces. Assume now that the matrix $M$ represents some inhomogeneous diffusion operator which contains information on the locations of the edges of the solution, cf. (1.13). The space $\mathcal{K}^{A^{\mathrm{T}}A}$ contains obviously no prior information, so it is solely determined through the action of $A^{\mathrm{T}}A$ on the right hand side $A^{\mathrm{T}}g$. In the space $\mathcal{K}^{A^{\mathrm{T}}A+\tau M}$, the prior information is superposed with the effect of applying $A^{\mathrm{T}}A$. Given that the features with highest a priori probabilities correspond to eigenvalues of $M$ that have the smallest magnitudes, a large number of iterations is required before $\mathcal{K}^{A^{\mathrm{T}}A+\tau M}$ is able to represent the edges. On the other hand, each of the Krylov vectors in $\mathcal{K}^{M^{-1}A^{\mathrm{T}}A}$ is directly influenced by the prior. As $M^{-1}$ amounts to steady-state solution operator of an inhomogeneous diffusion equation, each vector in the Krylov subspace $\mathcal{K}^{M^{-1}A^{\mathrm{T}}A}$, after an application of $A^{\mathrm{T}}A$, is subjected to inhomogeneous diffusion and so it contains the jump singularities used in the construction of $M$. Furthermore, notice that as the expected features are associated with the smallest magnitude eigenvalues of $M$, they correspond to the largest magnitude eigenvalues of $M^{-1}$.

In contrast to preconditioning, no clustering of eigenvalues takes place in priorconditioning. On the contrary, the condition number of the priorconditioned matrix can be larger than that of the original matrix. The working principle is fusing the information about the discontinuities into the forward operator so that the resulting Krylov subspace already in the initial iterations contains images of high prior probability with edges at the locations indicated by the prior. As a result, the priorconditioned algorithm usually requires only a few iterations to converge. While each priorconditioned iteration is more expensive than its unpriorconditioned counterpart due to the application of $M^{-1}$, for many large-scale problems the additional cost is insignificant and outweighted by the benefit of faster convergence, see the discussion in Section 6.

**2.2. Example problem: deconvolution in 1D.** Throughout the paper we use the following 1D deconvolution problem to exemplify the intermediate results. Assume that the matrix $A$ in the observation model

$$(2.4) \qquad g = Af + n$$

is a discretization of the stationary convolution operator

$$\mathcal{A}f(x) := \int_{-\infty}^{\infty} \mathcal{K}_{\sigma_f}(|x-x'|)f(x')\mathrm{d}x', \qquad \mathcal{K}_{\sigma_f}(x) = \sqrt{\frac{2}{\pi\sigma_f^2}}\exp\left[-\frac{x^2}{2\sigma_f^2}\right],$$

and $n \sim \mathcal{N}(0, \sigma_n^2 I)$ represents isotropic Gaussian white noise with zero mean. Figure 2.1(a) shows the target function $f$ along with its convoluted, noisy version $g$ with $\sigma_f = 0.03$, $\sigma_n = 0.01$. The functions were evaluated on a regular grid with 512 samples over the interval $[0, 1]$.

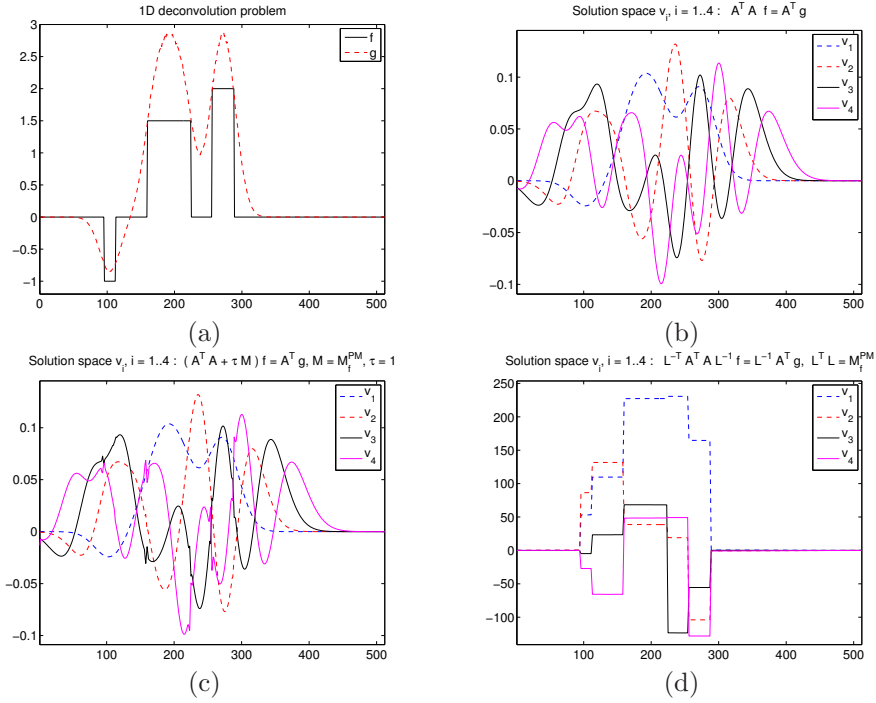

FIG. 2.1. *(a) true solution $f$ and data $g$. First four basis vectors of (b) $\mathcal{K}^{A^\mathrm{T} A}$ ($\tau = 0$), (c) $\mathcal{K}^{A^\mathrm{T} A + \tau M}$ ($\tau = 1$, applied to (1.8)), (d) $\mathcal{K}^{M^{-1} A^\mathrm{T} A}$ (applied to (1.7)).*

To illustrate the effect of priorconditioning claimed in Section 2.1, we solve (2.4) regularized with an (unrealistic) ideal regularizer $M$, discretization of $-\nabla \cdot c_{f_\mathrm{true}} \nabla$ c.f. (1.13), based on the *true* solution $f = f_\mathrm{true}$. Here, we used the Perona–Malik regularizer (1.5) (left equation) with $T = 0.005$. Practical cases will be discussed in Sections 4 and 5 where we show examples with fully nonlinear regularizers. The first four basis vectors of the respective Krylov spaces $\mathcal{K}^{A^\mathrm{T} A}$, $\mathcal{K}^{A^\mathrm{T} A + \tau M}$ (with $\tau = 1$) and $\mathcal{K}^{M^{-1} A^\mathrm{T} A}$ are shown in Figure 2.1(b–d). The vectors in $\mathcal{K}^{A^\mathrm{T} A}$ do not contain any prior information on the discontinuities, and so they are completely smooth. The third and fourth basis vectors in $\mathcal{K}^{A^\mathrm{T} A + \tau M}$ start to show small oscillations around the jumps of the true solution superposed with smooth components. The oscillations are, however, faint, which is consistent with the form of $A^\mathrm{T} A + \tau M$. As each vector in the Krylov subspace $\mathcal{K}^{M^{-1} A^\mathrm{T} A}$, after an application of $A^\mathrm{T} A$, is multiplied by $M^{-1}$, it contains the edges used in the construction of $M$.

**3. Implicitly preconditioned LSQR.** The split preconditioned symmetric system

$$(3.1) \qquad L^{-\mathrm{T}} A^\mathrm{T} A L^{-1} \hat{f} = L^{-\mathrm{T}} A^\mathrm{T} g, \qquad f = L^{-1} \hat{f}$$

can be solved by preconditioned conjugate gradient (PCG) without the need to provide a factorization $L^\mathrm{T} L$ for $M$ [42]. The key is to use suitable scalar products. For the left preconditioned normal equation

$$(3.2) \qquad M^{-1} A^\mathrm{T} A f = M^{-1} A^\mathrm{T} g,$$

PCG utilizes the $M$-weighted scalar product in which $M^{-1} A^\mathrm{T} A$ is self-adjoint:

$$\begin{aligned}(M^{-1} A^\mathrm{T} A f, g)_M &= (M^{-1} A^\mathrm{T} A f, M g) = (A^\mathrm{T} A f, g) \\ &= (f, A^\mathrm{T} A g) = (f, M M^{-1} A^\mathrm{T} A g) = (f, M^{-1} A^\mathrm{T} A g)_M.\end{aligned}$$

Alternatively, the right preconditioned normal equation

$$(3.3) \qquad A^\mathrm{T} A M^{-1} \bar{f} = A^\mathrm{T} g, \qquad M^{-1} \bar{f} = f,$$

can be solved by preserving symmetry with the $M^{-1}$-weighted scalar product, in which $A^\mathrm{T} A M^{-1}$ is self-adjoint. It is easy to see that both the left and right preconditioned variants produce the same sequence of computations and hence they are analytically equivalent.

**3.1. Preconditioned LSQR.** The LSQR algorithm [37, 36] is analytically equivalent to the conjugate gradient method applied to the normal equation, but it avoids explicit formation of the normal equation at any stage. Instead, LSQR makes use of the Golub–Kahan bidiagonalization [17]. Applied to the least squares problem

$$(3.4) \qquad \min_{f} \|g - Af\|$$

with a starting vector $g$, the bidiagonalization procedure can be written in matrix form as

$$\begin{aligned} U_{i+1}(\beta_1 e_1) &= g \\ A V_i &= U_{i+1} B_i \\ A^\mathrm{T} U_{i+1} &= V_i B_i^\mathrm{T} + \alpha_{i+1} v_{i+1} e_{i+1}^\mathrm{T}, \end{aligned}$$

where $e_i$ denotes the $i$th canonical basis vector, $\alpha_i \geq 0$ and $\beta_i \geq 0$ are chosen such that $\|u_i\| = \|v_i\| = 1$ and

$$(3.5) \qquad B_i = \begin{bmatrix} \alpha_1 & & & \\ \beta_2 & \alpha_2 & & \\ & \beta_3 & \ddots & \\ & & \ddots & \alpha_i \\ & & & \beta_{i+1} \end{bmatrix}, \qquad U_i = [u_1, u_2, \ldots, u_i], \quad V_i = [v_1, v_2, \ldots, v_i].$$

As all the quantities $B_i$, $U_i$ and $V_i$ are independent of $\tau$, the projected least squares problem

$$(3.6) \qquad \min_{y_i} \left\| \begin{bmatrix} B_i \\ \sqrt{\tau} I \end{bmatrix} y_i - \beta_1 e_1 \right\|,$$

is a proper generalization of $\min_{y_i} \|B_i y_i - \beta_1 e_1\|$ when $\tau = 0$. The projected problem (3.6) is then solved using QR decomposition yielding the approximation for the solution of the original problem, $f_i = V_i y_i$.

In the following, we derive a factorization-free priorconditioned LSQR algorithm, which is analytically equivalent to the PCG algorithm with the $M$-weighted inner products applied to the left preconditioned normal equation (3.2). The same algorithm (same sequence of computations) can be derived based on PCG with $M^{-1}$-weighted inner products applied to right preconditioned normal equation (3.3). We chose to present the $M$-weighted variant (left preconditioning) because it works with the original solution, unlike the $M^{-1}$-weighted one (right preconditioning) which involves a change of basis.

Recall that the preconditioned LSQR solves the problem where the preconditioner is applied to the least squares problem

$$\hat{f} = \operatorname{argmin} \|g - AL^{-1} \hat{f}\|,$$

resulting in the normal equation (3.1), which is exactly the split preconditioned normal equation (2.3) without damping. Hence, the preconditioned LSQR written out in Algorithm 1 is analytically equivalent to CG applied to the split preconditioned normal equation.

---
**Algorithm 1** Preconditioned LSQR, mathematically equivalent to CG applied to split preconditioned normal equation.
---
1: **Initialization:**
2: $\beta_1 u_1 = g$
3: $\alpha_1 v_1 = L^{-T} A^T u_1$
4: $w_1 = v_1$, $\hat{f}_0 = 0$, $\bar{\phi}_1 = \beta_1$, $\bar{\rho}_1 = \alpha_1$
5: **for** $i = 1, 2, \ldots$ **do**
6:    **Bidiagonalization:**
7:    $\beta_{i+1} u_{i+1} = A L^{-1} v_i - \alpha_i u_i$
8:    $\alpha_{i+1} v_{i+1} = L^{-T} A^T u_{i+1} - \beta_{i+1} v_i$
9:    **Orthogonal transformation:**
10:    $\rho_i = (\bar{\rho}_i^2 + \beta_{i+1}^2)^{1/2}$
11:    $c_i = \bar{\rho}_i / \rho_i$
12:    $s_i = \beta_{i+1} / \rho_i$
13:    $\theta_{i+1} = s_i \alpha_{i+1}$
14:    $\bar{\rho}_{i+1} = -c_i \alpha_{i+1}$
15:    $\phi_i = c_i \bar{\phi}_i$
16:    $\bar{\phi}_{i+1} = s_i \bar{\phi}_i$
17:    **Update:**
18:    $\hat{f}_i = \hat{f}_{i-1} + (\phi_i / \rho_i) w_i$
19:    $w_{i+1} = v_{i+1} - (\theta_{i+1} / \rho_i) w_i$
20:    **Break if stopping criterion satisfied**
21: **end for**
22: Transform back to original solution: $f = L^{-1} \hat{f}$
---

**3.2. Factorization-free preconditioning.** To derive the factorization-free priorconditioned variant of LSQR corresponding to the left preconditioned normal equation, we introduce new variables $v_i = L \tilde{v}_i$ and $w_i = L \tilde{w}_i$, cf. Algorithm 1. We then reformulate Algorithm 1 in terms of these new variables. The only parts affected are steps 3–4 of the initialization, the bidiagonalization and the update stage. Observing that step 3 of the initialization is of the same form (with $v_0 = 0$) as step 8 of the bidiagonalization, it is sufficient to consider only the latter. As $\hat{f}_i$ is a linear combination of $v_j$, $j \leq i$, the transformed solution has to be a linear combination of $\tilde{v}_j = L^{-1} v_j$, $j \leq i$. Thus, the resulting algorithm will directly produce a sequence of approximate solutions $f_i = L^{-1} \hat{f}_i$ and the change of variables at step 22 of Algorithm 1 cancels out.

The bidiagonalization in terms of $\tilde{v}_i$, $\tilde{w}_i$ and the $M$-weighted inner product reads

$$\beta_{i+1} u_{i+1} = A \tilde{v}_i - \alpha_i u_i$$
$$\tilde{v}_{i+1} = M^{-1} A^T u_{i+1} - \beta_{i+1} \tilde{v}_i$$
$$\alpha_{i+1} = ((\tilde{v}_{i+1}, \tilde{v}_{i+1})_M)^{1/2}$$
$$\tilde{v}_{i+1} = \tilde{v}_{i+1} / \alpha_{i+1}.$$

Keeping just one additional vector $\tilde{p} = M \tilde{v}_{i+1}$, the factorization-free algorithm needs no multiplications by the priorconditioner $M$, just one solve with $M$ per iteration. This amounts to rephrasing the bidiagonalization in the following way

$$\beta_{i+1} u_{i+1} = A \tilde{v}_i - \alpha_i u_i$$
$$\tilde{p} = A^T u_{i+1} - \beta_{i+1} \tilde{p}$$
$$\tilde{v}_{i+1} = M^{-1} \tilde{p}$$
$$\alpha_{i+1} = (\tilde{v}_{i+1}, \tilde{p})^{1/2}$$
$$\tilde{p} = \tilde{p} / \alpha_{i+1}$$
$$\tilde{v}_{i+1} = \tilde{v}_{i+1} / \alpha_{i+1}.$$

To conclude, the update in terms of the new variables and the original solution reads

$$f_i = f_{i-1} + (\phi_i / \rho_i) \tilde{w}_i$$

$$\tilde{w}_{i+1} = \tilde{v}_{i+1} - (\theta_{i+1}/\rho_i)\tilde{w}_i.$$

The resulting factorization-free preconditioned LSQR method (MLSQR) is summarized in Algorithm 2.

---

**Algorithm 2** MLSQR: Factorization-free preconditioned LSQR.

---
1: **Initialization:**
2: $\beta_1 u_1 = g$
3: $\tilde{p} = A^{\mathrm{T}} u_1$
4: $\tilde{v}_1 = M^{-1}\tilde{p}$
5: $\alpha_1 = (\tilde{v}_1, \tilde{p})^{1/2}$
6: $\tilde{v}_1 = \tilde{v}_1/\alpha_1$
7: $\tilde{w}_1 = \tilde{v}_1,\ f_0 = 0,\ \bar{\phi}_1 = \beta_1,\ \bar{\rho}_1 = \alpha_1$
8: **for** $i = 1, 2, \ldots$ **do**
9:    **Bidiagonalization:**
10:    $\beta_{i+1} u_{i+1} = A\tilde{v}_i - \alpha_i u_i$
11:    $\tilde{p} = A^{\mathrm{T}} u_{i+1} - \beta_{i+1}\tilde{p}$
12:    $\tilde{v}_{i+1} = M^{-1}\tilde{p}$
13:    $\alpha_{i+1} = (\tilde{v}_{i+1}, \tilde{p})^{1/2}$
14:    $\tilde{p} = \tilde{p}/\alpha_{i+1}$
15:    $\tilde{v}_{i+1} = \tilde{v}_{i+1}/\alpha_{i+1}$
16:    **Orthogonal transformation:**
17:    Steps 10-16 as in Algorithm 1
18:    **Update:**
19:    $f_i = f_{i-1} + (\phi_i/\rho_i)\tilde{w}_i$
20:    $\tilde{w}_{i+1} = \tilde{v}_{i+1} - (\theta_{i+1}/\rho_i)\tilde{w}_i$
21:    **Break if stopping criterion satisfied**
22: **end for**

---

**3.3. LSQR with regularization.** Two types of regularization are relevant in our framework: Tikhonov regularization, where the parameter $\tau$ controls the amount of regularization, and early truncation of Krylov methods, in which regularization arises from the problem being projected into a small dimensional subspace.

When priorconditioning is used, Tikhonov regularization results in a simple damped least squares problem (2.3). Damping for a fixed value of $\tau$ can be easily incorporated in LSQR as described in [36] at the cost of doubling the number of Givens rotations. Due to the shift invariance of Krylov spaces, different choices of the damping parameter $\tau$ result in the same Krylov subspace and Lanczos vectors $V_i$. If $V_i$ are stored, the projected least squares problem (3.6) can be efficiently solved for multiple values of $\tau$, which is of benefit when the value of $\tau$ is not known in advance. However, the additional storage requirements may limit the feasibility of such an approach. Some techniques for solving (3.6) with a variable $\tau$ are discussed in [4] using singular values and the first and last rows of the matrix of the right singular vectors of the bidiagonal matrix $B_i$. Those quantities can be obtained at the cost $\mathcal{O}(i^2)$ at the $i^{\text{th}}$ iteration. These strategies are however only viable for a limited number of iterations $i$ as the singular value decomposition of the bidiagonal matrix $B_i$ can not be efficiently updated even though $B_i$ simply expands by a row and a column in each iteration. For larger $i$, the algorithm described in [11] for the least squares solution of (3.6) at the cost of $\mathcal{O}(i)$ for each value of $\tau$ is the preferable option.

**3.4. Stopping criteria.** The original paper on LSQR [37] discusses three stopping criteria:
   S1: $\|\bar{r}_i\| \leq \text{BTOL}\|g\| + \text{ATOL}\|\bar{A}\|\|f_i\|$ (consistent systems),
   S2: $\frac{\|\bar{A}^T \bar{r}_i\|}{\|\bar{A}\|\|\bar{r}_i\|} \leq \text{ATOL}$ (inconsistent systems),
   S3: $\text{cond}(\bar{A}) \geq \text{CONLIM}$ (both),

where $\bar{r}_i := \bar{g} - \bar{A}f_i$ is the residual of (1.8) with

$$\bar{A} = \begin{bmatrix} A \\ \sqrt{\tau}L \end{bmatrix}, \qquad \bar{g} = \begin{bmatrix} g \\ 0 \end{bmatrix},$$

and ATOL, BTOL and CONLIM are user-specifiable parameters, refer to [37] for details.

It is well known that LQSR applied to the least squares problem (1.8) monotonically decreases the residual norm $\|\bar{r}_i\|$ since LSQR is analytically equivalent to the CG method applied to the corresponding normal equation (2.1). For ill-posed problems, we also consider the Morozov discrepancy principle

S4: $\|r_i\| \leq \eta \delta$,

where $r_i := g - Af_i$ and $\delta$ is the (estimated) noise level, see e.g. [22] for related considerations on regularization properties of CG type methods. The factor $\eta > 1$ is included to prevent underregularization [9]. Whenever $\tau = 0$, it holds that $r_i = \bar{r}_i$ and the sequence $\|r_i\|$ is also monotonically decreasing.

While the stopping criterion S4 is not affected by the priorconditioning (the residual $r$ remains the same) as priorconditioning effectively substitutes $\bar{A}L^{-1}$ for $\bar{A}$, the criteria S1–S3 monitor different values for the same solution $f$ depending on whether $f$ was obtained with priorconditioning or not.

**3.5. MLSQR: example 1D deconvolution problem.** To demonstrate the difference in solutions obtained with and without priorconditioning, we revisit the linear deconvolution problem of Section 2.2 still using the ideal linear regularizer resulting in preconditioner $M := M_{f_{\text{true}}} = L^T L$.

We solve the least squares problem (1.8) with and without priorconditioning using the factorization-free LSQR. The regularization parameter value choice of $\tau = 1$ was made by inspection and was optimized for the unpriorconditioned problem. In addition, the value $\tau = 0$ was tested. For the priorconditioned problem selecting $\tau = 0$ results in a priorconditioned solution without damping, while for the unpriorconditioned problem setting $\tau = 0$ eliminates (Tikhonov) regularization.

Figures 3.1(a,b) show the solutions obtained with and without priorconditioning. In Figure 3.1(a), we selected the solutions using stopping criterion S4 with $\delta = 10^{-2}$ and $\eta = 1.1$, which resulted in stopping after 9 iterations for the priorconditioned problem and after 16 iterations for the unpriorconditioned problem for both tested values of $\tau$. The discrepancy principle S4 provides a good stopping criterion for the priorconditioned problem, but terminates the LSQR algorithm for the regularized unpriorconditioned problem too early. Figure 3.1(b) shows the solutions obtained with the stopping criterion S2 with ATOL $= 10^{-2}$. While S2 is useful for stopping the unpriorconditioned regularized problem, it stops the priorconditioned problem prematurely with the chosen ATOL. Comparing the best respective solutions for both formulations, the priorconditioned formulation after 9 iterations yields a better solution than the unpriorconditioned regularized formulation after 32 iterations.

**4. Least squares with nonlinear regularizer.** When employing a general regularizer of the form (1.2), the computation of the MAP estimate becomes a nonlinear least squares optimization problem (1.6). In this section, we propose a method for solving (1.6) based on lagged diffusivity fixed point iteration and the factorization-free priorconditioned LSQR algorithm.

**4.1. Nonlinear solver.** The necessary condition for the minimizer of (1.6), $\nabla \Phi = 0$, results in the nonlinear equation

$$(4.1) \qquad (A^T A + \tau M_f) f = A^T g,$$

where $M_f$ is an approprioate discretization of (1.13). We solve (4.1) using the lagged diffusivity method [52]

$$(4.2) \qquad (A^T A + \tau M_{f^k}) f^{k+1} = A^T g,$$

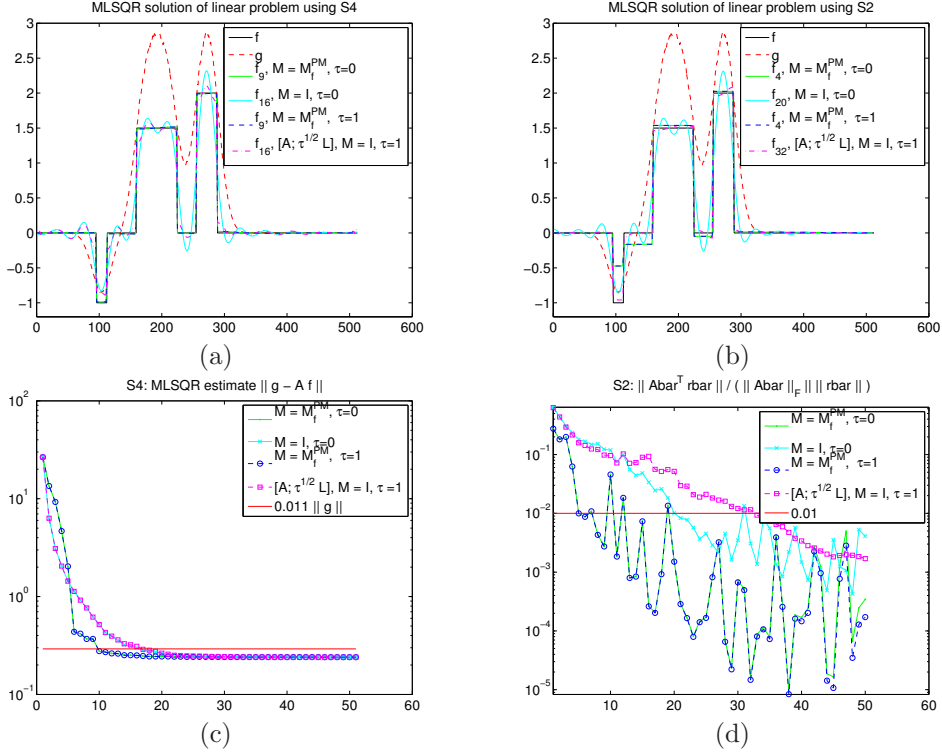

Fig. 3.1. *Upper row: Solution of the linear deconvolution problem computed with LSQR with and without preconditioning (a) using stopping criterion S4, (b) using stopping criterion S2. Lower row: Behavior of the quantities used for stopping criteria. (c) S4: $\|g - Af\| \leq 1.1 \cdot 10^{-2}$. (d) S2: $\|\bar{A}^T \bar{r}\|/(\|\bar{A}\|\|\bar{r}\|) \leq 10^{-2}$.*

a fixed-point iteration that was originally devised for total variation image denoising.

The linear problem at the $k^{\text{th}}$ step of the fixed point iteration (4.2) can be solved efficiently and matrix-free using factorization-free priorconditioned LSQR. Algorithm 3 gives the details of the resulting method for solving the nonlinear least squares problem (1.6). Notice that there are two iteration levels: the *outer* arising from lagged diffusivity, and the *inner* corresponding to MLSQR.

---
**Algorithm 3** Priorconditioned lagged diffusion solver.
---
1: Initialize with $f^0 = 0$
2: **for** $i = 1, 2, \ldots$ **do**
3:     Form the preconditioner $M_{f^{k-1}}$ using the current approximation $f^{k-1}$
4:     Solve symmetrically preconditioned linear system (2.3)

$$(4.3) \qquad (L^{-T} A^T A L^{-1} + \tau I)\hat{f}^k = L^{-T} A^T g, \quad f^k = L^{-1} \hat{f}^k,$$

    where $M_{f^{k-1}} = L^T L$, with MLSQR and stopping criterion S4.
5:     **Break** if $R(f^k)$ does not decrease fast enough
6: **end for**
---

The objective of Algorithm 3 is to minimize the functional $\Phi$ in (1.6). Under certain assumptions and for a parameter value $\tau$ chosen such that the Morozov discrepancy principle S4 is satisfied, the minimization problem can be reformulated as the residual method [20],

$$(4.4) \qquad \min_{f \in X} R(f) \text{ subject to } \|g - Af\| \leq \eta \delta.$$

Motivated by this form, we choose to monitor the value of the penalty $R$ and stop the lagged diffusivity fixed point iteration when $R$ does not decrease fast enough. Because each iterate $f^k$ has approximately the same residual norm $\|r^k\| \approx \delta$, the value of $\Phi$ will in practice decrease until a feasible solution to (1.6) is found. We emphasize that it is important to stop the Krylov method in step 4 of Algorithm 3 consistently: We want to avoid the situation where the trend in values of $R(f^k)$ becomes distorted by the Krylov solver. (For example, taking few Krylov iterations in one step while taking many iterations in the next will likely lead to increase in $R(f^k)$.) In practice, it may take numerous Krylov iterations before the first few lagged diffusivity iterates reach the noise level. This issue is effectively circumvented by combining the Morozov discrepancy principle with a global upper limit on the number of allowed Krylov iterations within each lagged diffusivity step. Although the noise level is not reached in earlier iterations, it does not seem to affect the Algorithm 3, see Section 5.3.

A few more comments are in order regarding Algorithm 3:
- We deliberately choose to solve for $f^{k+1}$ rather than the update $\Delta f^{k+1} := f^{k+1} - f^k$. This is because the employed regularizers encode the information on the solution $f$ and not the updates $\Delta f^{k+1}$.
- Our simulations suggest that solutions of high quality can be obtained by setting $\tau = 0$ in (2.3) and employing early stopping of the MLSQR algorithm. A clear advantage of this approach is that the data discrepancy $\|r^k\|$ is being minimized instead of $\|\bar{r}^k\|$. However, setting $\tau = 0$ severs the connection between Algorithm 3 and the original problem of finding the MAP estimate for (1.1).
- If $M_f^{-1}$ is computed exactly, $\tau \neq 0$ is used and MLSQR is iterated until convergence, our scheme is equivalent to lagged diffusivity iteration. In practice, one or more of these conditions may not hold.

**4.2. Inverting the preconditioner.** While our preconditioner is nonstandard in the sense that it does not approximate the inverse of the forward operator, it still amounts to a solution of an elliptic partial differential equation. Therefore, after appropriate discretization, applying the priorconditioner is reduced to a problem of inverting a sparse matrix. There exists a bulk of literature discussing these topics, see e.g. [5, 29, 48] for treatises on partial differential equations and [42, 50] for methods applicable to solving linear problems in general.

In many applications, in particular in one or two dimensions, the preconditioner can simply be inverted through its Cholesky factorization in conjunction with a suitable reordering algorithm that keeps the fill-in small. For large problems, more sophisticated sparse direct solvers have been developed. A comparison of numerical packages for direct solution of sparse symmetric systems is given in [19]. These methods use graph models to decrease the computational load arising from the factorization and the subsequent solution step. However, for large three-dimensional problems, it may still be necessary to work out of core, i.e. using main storage. Many of the sparse direct methods offer an out-of-core option, but the performance is obviously affected by the main storage access.

Multigrid methods are widely regarded as the state of the art for elliptic problems [5]. It is well known that the basic iterative relaxation methods, including the Jacobi method and symmetric successive overrelaxation (SSOR), quickly reduce the high-frequency part of errors in the approximate solutions. However, low frequency errors are diminished only slowly. Multigrid methods address this issue by restricting the residual to a coarser grid, where the low frequency errors from a fine grid appear as high frequencies and, hence, are quickly diminished by relaxation methods. The obtained correction can then be interpolated into the finer grid, resulting in fast convergence also for the low frequency errors. Naturally, a multi-level hierarchy of grids can be used, hence the name (geometric) multigrid.

As an alternative to a geometric construction of a hierarchical set of grids, it is possible to proceed in a purely algebraic manner. Starting from the finest level, coarser levels are constructed on the basis of the strength of the connections between the nodes. In the absence of the underlying grid, the equivalent of the high-frequency error is the error that is efficiently eliminated by relaxation methods.

In particular, the use of Jacobi, SSOR or incomplete LU decomposition (ILU) smoothers results in a symmetric preconditioner if the same number of iterations is taken in the pre- and post-smoothing stages. By using a fixed number of smoother steps and grid-to-grid operations, one obtains a fixed approximation of a matrix inverse that can be used to priorcondition Krylov methods such as the factorization-free priorconditioned LSQR derived in Section 3. While using a flexible Krylov subspace method [18, 31, 41, 44, 46] for solving (4.3) in Algorithm 3 would remove the requirement of using multigrid with fixed number of steps, such methods do not arguably provide any benefits over the approach taken in this article. This is because, remarkably, even a low-cost single V-cycle multigrid approximation to $M^{-1}$ is enough to provide the benefits of priorconditioning, cf. Section 5, and consequently there is no need for adapting the number of cycles.

**4.3. Nonlinear solver: example 1D deconvolution problem.** We are now in a position to consider the nonlinear 1D deconvolution problem with nonlinear Perona–Malik regularization (1.5) (left equation). We solve the associated minimization problem (1.6) with Algorithm 3.

In this particular example, we chose to use the undamped solution, i.e., $\tau = 0$. We employed the discrepancy principle (S4) with the noise level $\delta = 10^{-2}$ and $\eta = 1.1$ as a stopping criterion for MLSQR. We stopped the lagged diffusivity fixed point iteration when the relative change in the functional $R$ dropped below 15%. As the considered problem is small, we used Cholesky factorization for inverting the preconditioner.

Figure 4.1 shows the evolution of the first six basis vectors of $\mathcal{K}^{M^{-1}A^{\mathrm{T}}A}$ for $M = M_{f^k}$, $k \in \{1, 4, 7, 10\}$. For $k = 1$, the priorconditioner is noninformative (homogeneous). Hence, all of the shown basis vectors are smooth and therefore replicate the jumps of the target function poorly. However, the priorconditioner quickly adapts to the edges emerging in the intermediate solutions $f^k$ as can be seen from the basis vectors for $k > 1$.

The solution at every third and the error norm at every lagged diffusivity step are illustrated in Figures 4.2(a) and 4.2(b). The residual norm over the inner and outer iterations is shown in Figure 4.3(a). Notice that every lagged diffusivity iteration resets the residual norm to the value $\|g\|$. This is due to the initialization of MLSQR with $f = 0$ at the beginning of every outer iteration to allow the solution to develop features compatible with the new priorconditioner. The updated priorconditioner allows faster reduction of the residual norm and produces an improved solution in fewer MLSQR iterations as shown in Figure 4.3(c).

Figure 4.3(b) presents the values obtained by the functional $R$ over the inner and outer iterations. We observe that, within the computation of any given lagged diffusivity iterate, the general trend for $R$ is to increase as MLSQR progresses until a saturation level is reached. On the other hand, the saturation level decreases as the outer iteration advances. At a certain point, the saturation level stops decreasing (or is not decreasing fast enough), see Figure 4.3(d), which is when Algorithm 3 is terminated. The plot of the $L_2$ norm of the solution error in Figure 4.2(d) demonstrates that the last step taken before the stopping criterion was triggered did not perceivably improve the solution, which corroborates our choice of the stopping criterion.

**5. A 3D image reconstruction problem.** In this section we apply our method to a large-scale 3D image reconstruction problem on an unstructured grid.

**5.1. Fluorescence Diffuse Optical Tomography.** We consider fluorescence Diffuse Optical Tomography (fDOT), which indirectly measures the conversion of strongly scattered (i.e. diffusively propagating) light from an excitation wavelength to a (longer) emission wavelength in the presence of fluorescent markers accumulating in regions of interest. The goal is to monitor cellular and subcellular functional activity. Although fDOT is a method mostly used for small animal research [32, 33, 53, 30], it is a promising technique with several medical applications such as detection, diagnosis and monitoring of human neoplasms, in particular breast tumors [10, 35, 49]. Due to the diffusive nature of light propagation in biological tissue, the image reconstruction

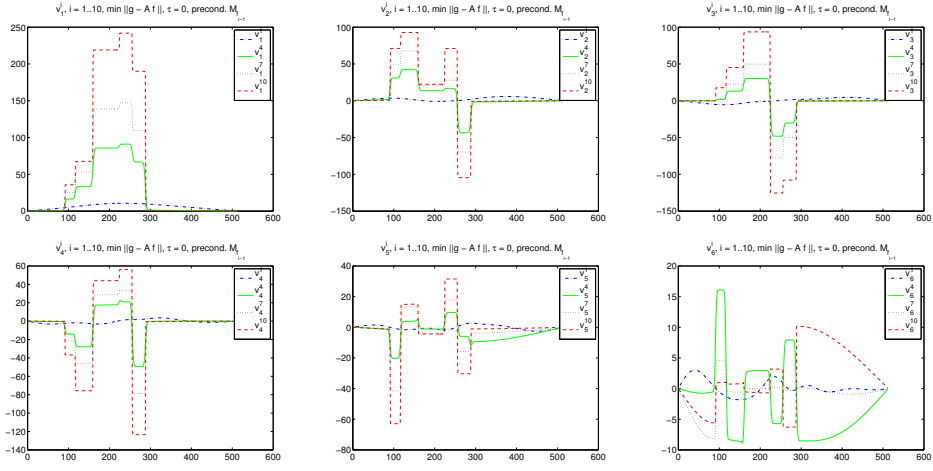

Fig. 4.1. *Evolution of the first six basis vectors of $\mathcal{K}^{M^{-1}A^{\mathrm{T}}A}$ through the outer iterations.*

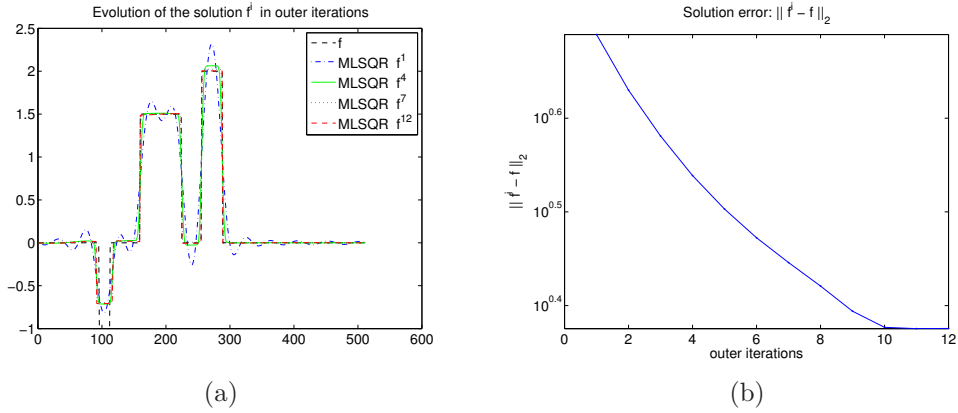

Fig. 4.2. *(a) Evolution of the solution computed with MLSQR through the outer iterations. (b) The $L_2$ norm of the solution error over outer iterations.*

is an ill-posed inverse problem [1, 16, 15].

In the linear approximation to fDOT, the forward model is described with coupled diffusion equations at the excitation and emission wavelengths, $\lambda_{\mathrm{e}}$ and $\lambda_{\mathrm{f}}$, respectively:

(5.1) $\quad [-\nabla \cdot \kappa(x, \lambda_{\mathrm{e}})\nabla + \mu_a(x, \lambda_{\mathrm{e}})] U(x, \lambda_{\mathrm{e}}) = 0, \quad x \in \Omega,$

(5.2) $\quad [-\nabla \cdot \kappa(x, \lambda_{\mathrm{f}})\nabla + \mu_a(x, \lambda_{\mathrm{f}})] U(x, \lambda_{\mathrm{f}}) = U(x, \lambda_{\mathrm{e}})h(x, \lambda_{\mathrm{f}}), \quad x \in \Omega.$

with $\mu_a$ and $\mu'_s$ the absorption and the reduced scattering coefficients, respectively, $\kappa = [3\,(\mu'_s + \mu_a)]^{-1}$ the diffusion coefficient and $h$ the fluorescence yield coefficient.

The diffusion equations are complemented by the respective Robin boundary conditions

(5.3) $\quad U(x, \lambda_{\mathrm{e}}) + 2\zeta\kappa(x, \lambda_{\mathrm{e}})\dfrac{\partial U(x, \lambda_{\mathrm{e}})}{\partial \nu(x)} = \Theta_{\mathrm{s}}(x)q, \quad x \in \partial\Omega,$

(5.4) $\quad U(x, \lambda_{\mathrm{f}}) + 2\zeta\kappa(x, \lambda_{\mathrm{f}})\dfrac{\partial U(x, \lambda_{\mathrm{f}})}{\partial \nu(x)} = 0, \qquad x \in \partial\Omega,$

where $\nu$ is the outward unit normal and $\zeta$ accounts for the refractive index mismatch at the boundary. The right-hand side in (5.3) models the effect of the excitation light source as an inward (diffuse) photon current, a product of the source emitted photon current $q$ and a source coupling coefficient function $\Theta_{\mathrm{s}}$ (notice no source term in (5.1)). On the other hand, the emission photon density arises solely from the fluorecence in $\Omega$ (right-hand side of (5.2)), resulting in homogeneous Robin boundary condition (5.4) .

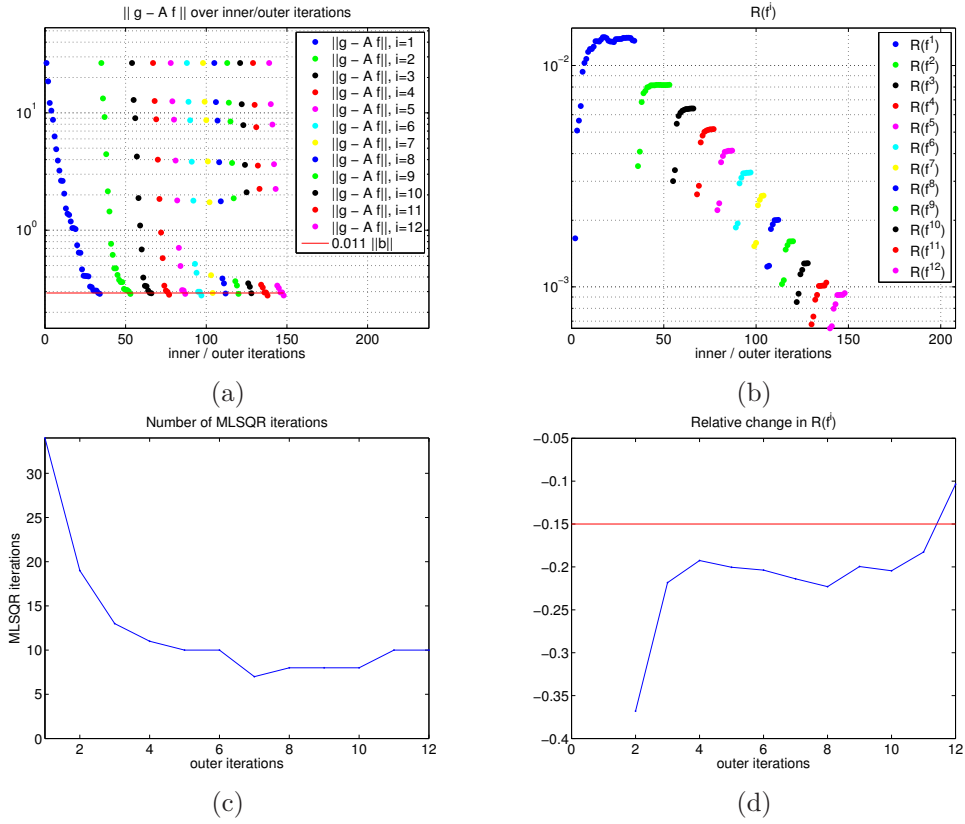

Fig. 4.3. *(a) Residual norm (S4) over lagged diffusivity/MLSQR iterations. (b) $R(f_i^k)$ over lagged diffusivity (k)/MLSQR (i) iterations. (c) Number of MLSQR iterations over lagged diffusion iterations. (d) Stopping criterion for the nonlinear iteration: $(R(f^k) - R(f^{k-1}))/R(f^{k-1})$ over lagged diffusivity iterations. The threshold level was chosen to be $-0.15$.*

The measured photon density for a detector with wavelength independent detector coupling coefficient $\Theta_d$ is given as

$$(5.5) \qquad y(\lambda) = \int_{\partial \Omega} \Theta_d(x) \frac{U(x,\lambda)}{2\zeta} dx, \quad \lambda \in \{\lambda_e, \lambda_f\}.$$

In practice several sources and detectors are deployed, resulting in a vector of measurements $y(\lambda) \in \mathbb{R}^{N_s N_d}$, where $N_s$ and $N_d$ is the number of sources and detectors, respectively. We furthermore use the following normalization, demonstrated to reduce the effects of unknown $\Theta_s$ and $\Theta_d$ [34, 45]

$$(5.6) \qquad g_j = \frac{y_j(\lambda_f)}{y_j(\lambda_e)}.$$

We now define the mapping $\mathcal{A} : h \mapsto g$ through solution of the system (5.1)–(5.6). $\mathcal{A}$ is clearly nonlinear in the parameters $\{\mu_a(x,\lambda_e), \mu_a(x,\lambda_f), \mu'_s(x,\lambda_e), \mu'_s(x,\lambda_f)\}$, but it is linear in the sought-for parameter $h(x,\lambda_f)$ (assuming that the other optical parameters are independent of $h$, which is reasonable for typical concentrations of fluorophores found in biomedical applications).

The adjoint mapping $\mathcal{A}^* : b \mapsto z$ for a single source-detector pair, is defined through solution of (5.1) with the boundary condition (5.3) and

$$(5.7) \qquad [-\nabla \cdot \kappa(x,\lambda_f)\nabla + \mu_a(x,\lambda_f)] U^*(x,\lambda_f) = 0, \quad x \in \Omega,$$

with the inhomogeneous Robin boundary condition

$$(5.8) \qquad U^*(x,\lambda_f) + 2\zeta \kappa(x,\lambda_f) \frac{\partial U^*(x,\lambda_f)}{\partial \nu(x)} = \frac{\Theta_d(x) b}{y(\lambda_e)}, \quad x \in \partial\Omega,$$

followed by the multiplication of the two solutions, yielding

$$(5.9) \qquad z = U^*(x, \lambda_\text{f}) U(x, \lambda_\text{e}).$$

If multiple detectors are used, the right-hand side of (5.8) will involve summation over the detectors. Analogously, a multiple source configuration will involve $N_\text{s}$ solves for (5.1), (5.3) and (5.7), (5.8), with (5.9) summing the pairwise products of $U$ and $U^*$ over the sources.

**5.2. Simulation setup.** Our test phantom is a cylinder of radius 25 mm and height 50 mm. The fluorophore distribution is represented through three spherical inclusions, one with $h = 0.06$ and the other two with $h = 0.1$, and the background value $h = 0$, see Figures 5.1, 5.2(e). The remaining optical parameters are assumed homogeneous, and their values are summarized in Table 5.1. The fluorophore was excited with each of 80 sources uniformly distributed along five rings on the boundary of the cylinder, see Figure 5.1. The measurements were sampled by 80 detectors placed along the same rings on the boundary, half way in between the sources, resulting in a total of 6400 measurements. All the involved coupling coefficient functions were modeled as Gaussian distributions on the boundary $\partial \Omega$. Measurement errors were simulated using additive Gaussian noise with standard deviation equal to 1% of the corresponding ideal measurement.

| | | |
|---|---|---|
| absorption coefficient $\mu_a(\cdot, \lambda_\text{e}) = \mu_a(\cdot, \lambda_\text{f})$ | 0.05 | mm$^{-1}$ |
| reduced scattering coefficient $\mu'_s(\cdot, \lambda_\text{e}) = \mu'_s(\cdot, \lambda_\text{f})$ | 1 | mm$^{-1}$ |
| refractive index | 1.4 | |

TABLE 5.1
*Optical parameters used in the simulations.*

The fDOT problem was modeled and discretized using the finite element method (FEM) with software package TOAST [2]. The discretization with first order Lagrangian elements on a tetrahedral mesh yielded 27084 degrees of freedom which corresponds to approximately 1.5mm resolution. The associated inverse problem was regularized with the differentiable approximation to the total variation functional in (1.2). The resulting matrix $M_f$ is thus a discrete representation of the operator $\nabla \cdot (|\nabla f|^2 + T^2)^{-1/2} \nabla$. We chose $T = 10^{-6}$ as it yields a well-conditioned matrix $M$ for fluorescence yield coefficients with edges of height approximately equal or smaller than the expected maximum of 0.1. The preconditioner was applied using a computationally low-cost algebraic multigrid solver with single V-cycle and two steps of ILU(0) smoothing, implemented in IFISS [43, 13].

**5.3. Solution of the fDOT problem.** We solve the fDOT problem described in Section 5.1 using Algorithm 3. The proposed method is particularly well suited for problems of this kind: fDOT is large-scale and the explicit construction of the matrix $A$ representing the discretized forward mapping can be impractical. Instead, the action of $A$ on a vector is obtained through solution of (5.1)–(5.6) for each source. Similarly, multiplication by $A^\text{T}$ involves solving (5.1), (5.3) and (5.7)–(5.9). In our example, both of these mappings amount to 160 solves of elliptic partial differential equations. Moreover, because the FEM mesh is unstructured, the priorconditioner $M$ naturally arises in the unfactorized form.

We demonstrate the benefits of priorconditioning by comparing the performance of Algorithm 3 to an unpriorconditioned reference method which is an adaptation of Algorithm 3 obtained by substituting the unpriorconditioned normal equation (4.2) for the priorconditioned one (4.3) at step 4. Since (4.2) is well-posed for appropriately chosen $\tau$ and the corresponding linear system inconsistent, the unpriorconditioned variant is stopped with S2 with ATOL = $10^{-3}$, while the priorconditioned variant with S4 and $\eta = 1.1$. Both methods use the regularization parameter value $\tau = 10^4$, which was tuned to yield the best-case results for the unpriorconditioned reference method. As in Section 4.3, we also test $\tau = 0$ and show that the proposed algorithm is very robust with respect to choice of $\tau$. Note that for $\tau = 0$ the unpriorconditioned variant is reduced to a linear least squares problem (3.4) solved by LSQR. Since (3.4) is ill-posed, in this case LSQR is stopped using the criterion S4 (as in the priorconditioned cases).

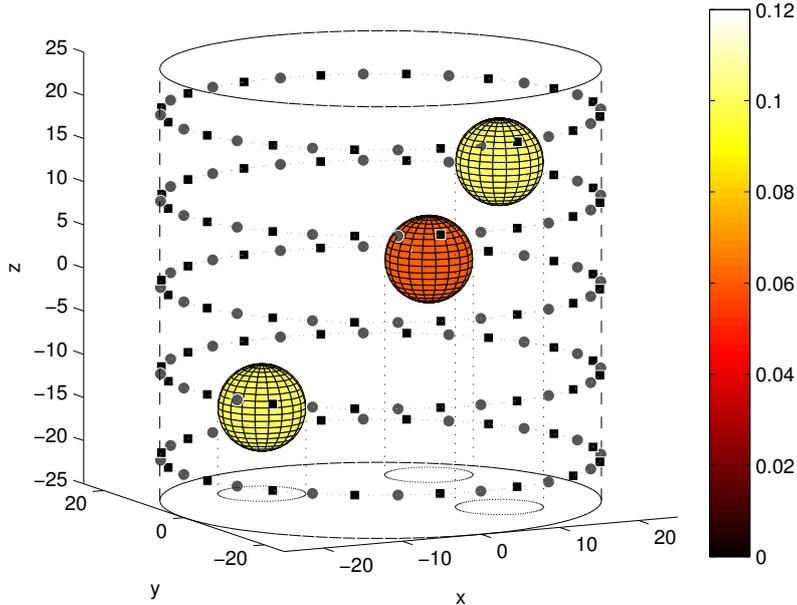

Fig. 5.1. *The cylindrical phantom with three inclusions. Source and detector locations are marked with gray points and black squares, respectively.*

We monitor the value of the penalty $R(f^k)$ of the intermediate solutions $f^k$ and terminate the lagged diffusivity (outer) iteration in Algorithm 3 when $R(f^k)$ stops decreasing. Such stopping criterion is not viable for the unpriorconditioned variant because the value of the penalty evolves differently for the two methods. Therefore, a fixed number of 25 lagged diffusivity steps was taken to ensure convergence for the unpriorconditioned algorithm with $\tau = 10^4$. The resulting reconstructions are depicted in Figures 5.2 and 5.3. The two unpreconditioned solutions suffer from overshooting and spurious oscillations, while their priorconditioned counterparts have the correct shape and estimate the value of the fluorescence yield coefficient $h$ more accurately. These findings are also reflected by the $L_2$ norm of the error of the solutions: The unpriorconditioned algorithm attains error norm of 0.7158 and 1.5025 for $\tau = 10^4$ and $\tau = 0$, respectively. The corresponding error norm values for the priorconditioned variant were 0.5428 and 0.5407. The error norms in each lagged diffusivity iteration are plotted in Figure 5.4(a).

Figure 5.4(b) shows the number of MLSQR iterations taken within each lagged diffusivity step. We chose to limit the number of the inner iterations for the priorconditioned algorithm to a maximum of 20 in accordance with the discussion in Section 4.1. The maximal allowed number of inner iterations was attained in the first 6 ($\tau = 10^4$) and 5 ($\tau = 0$) lagged diffusivity iterations, while MLSQR was stopped by the Morozov criterion in fewer than 20 iterations in the following lagged diffusivity steps. In fact, the motivation for limiting the maximum number of MLSQR iterations is that the first one or two lagged diffusivity steps would need an disproportionate number of MLSQR iterations to reach the residual level suggested by the Morozov discrepancy principle. The limit on Krylov iterations alleviates this issue without a noticeable effect on the final reconstruction. While after the initial phase the number of MLSQR iterations decreases to 9 for the priorconditioned algorithm, the unpriorconditioned variant requires an increasing number of Krylov steps, which we attribute to the growing condition number of the regularizer $M_{f^k}$ over the lagged diffusivity iterations.

Figure 5.4(c) shows the evolution of the penalty term $R(f)$ over lagged diffusivity iterations. In the first outer iteration, the two priorconditioned solutions show a pronounced difference. The initial priorconditioner is a discretization of the Laplacian and contains no information about the edges in the solution. For this reason, the damping flattens the first outer iterate particularly strongly, and consequently the associated penalty $R(f) = \int_\Omega |\nabla f(x)| \mathrm{d}x$ is of smaller magnitude than it is for the undamped variant. Notice that similar flattening takes place for the unpriorconditioned method, as well. Once the priorconditioner contains some information about the edges, the behavior of both priorconditioned cases is qualitatively the same. Apart from the first lagged diffusivity iteration, the penalty decreases over the lagged diffusivity iterations, providing further evidence to support the penalty-based stopping criterion of Algorithm 3.

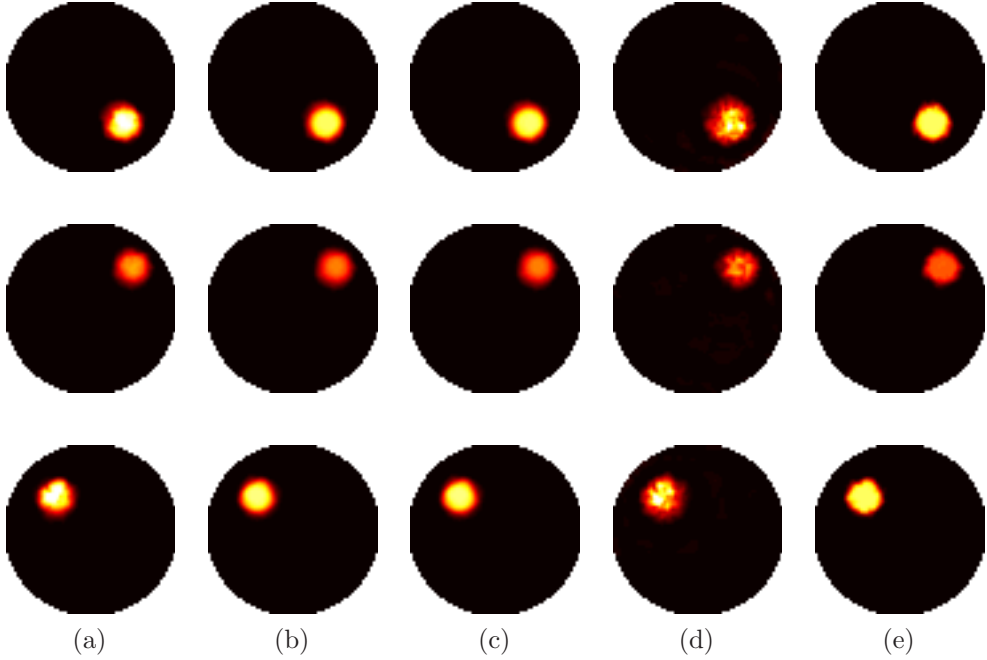

(a) (b) (c) (d) (e)

FIG. 5.2. *Transaxial planes through the centers of the spherical inclusions at $z = 15$, $z = 0$ and $z = -15$ of the solution to (a) unpriorconditioned problem, (b) priorconditioned problem with damping $\tau = 10^4$, (c) priorconditioned problem without damping $\tau = 0$, and (d) problem with no other regularization but early stopping of LSQR; (e) the phantom. (x-axis horizontal, y-axis vertical).*

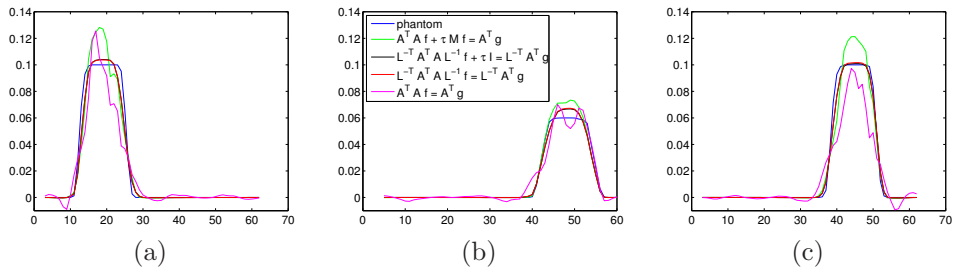

(a) (b) (c)

FIG. 5.3. *Cross sections along the x-axis through the centers of the inclusions of the solutions in Figure 5.2 (a) $z = -15, y = 10$, (b) $z = 0, y = 12.5$ and (c) $z = 15, y = -10$.*

**6. Conclusions.** In this paper, we considered efficient solution of large-scale linear ill-posed inverse problems with nonlinear regularization. We devised a highly efficient matrix-free algorithm for solution of such problems combining a lagged diffusivity fixed point iteration with priorconditioning of Krylov methods. Priorconditioning affords a way of embedding the information contained in the prior directly into the forward operator resulting in highly accelerated convergence of Krylov methods. A novel factorization-free preconditioned LSQR algorithm was presented for solving the linear priorconditioned problem which allows an implicit application of the preconditioner through efficient schemes such as multigrid. This is of particular interest for problems formulated on unstructured grids, where the preconditioner naturally occurs in an unfactorized form and the factorization is computationally infeasible. Furthermore, the relevant regularizers arise as discretizations of elliptic partial differential equations for which approaches like multigrid have been extensively studied and applied.

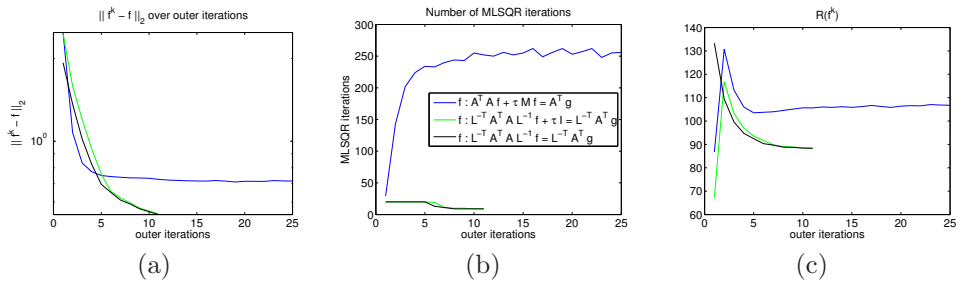

Fig. 5.4. *(a) Error norm plotted over outer iterations. (b) Number of MLSQR iterations plotted over outer iterations. (c) Penalty R plotted over outer iterations. (The unpriorconditioned method with $\tau = 0$ does not involve outer iterations. Therefore, the associated results are omitted from this figure.)*

The presented algorithm is matrix-free, i.e. capable of solving problems where the forward mapping cannot be computed and/or stored explicitly as a matrix. In particular, in such cases the cost of the application of the forward and adjoint mappings in each step of the Krylov method overbears the additional cost of preconditioning, while the use of priorconditioning reduces the number of Krylov iterations to a fraction of those needed without priorconditioning.

We demonstrated the effectiveness of our algorithm on the 3D fDOT image reconstruction problem. For truly large-scale problems, for which the diffusion equations involved in the forward and adjoint mapping have to be solved iteratively, the cost of applying the forward and adjoint mappings is typically running hundreds of times higher than the cost of applying the preconditioner, making the reduction in the number of Krylov iterations paramount.

The present paper deals with ill-posed problems, where the forward mapping is linear and the nonlinearity is limited to the regularization term. We intend to investigate the extension of the ideas to fully nonlinear inverse problems such as electrical impedance tomography and diffuse optical tomography.

**7. Acknowledgements.** The authors would like to acknowledge the financial support of the funding bodies: Marta Betcke is supported by an EPSRC Postdoctoral Fellowship (grant number EP/H02865X/1) and Lauri Harhanen by Finnish Doctoral Programme in Computational Sciences (FICS) and Finnish Doctoral Program in Inverse Problems.